\documentclass[leqno, 10pt]{amsart}
\usepackage{amsmath}
\usepackage{amssymb, latexsym}
  \newtheorem{definition}{Definition}[section]
 \newtheorem{theorem}{Theorem}[section]
 
 \newtheorem{proposition}{Proposition}[section]

 \newtheorem*{theorem*}{Theorem}
 
 \newtheorem*{proposition*}{Proposition}

 \theoremstyle{remark}

\newcommand{\op}[1]{\operatorname{#1}}

\newcommand{\brak}[1]{\ensuremath{\langle\! #1\!\rangle}}

\newcommand{\C}{\ensuremath{\mathbb{C}}} 
 
\newcommand{\bH}{\ensuremath{\mathbb{H}}} 
\newcommand{\R}{\ensuremath{\mathbb{R}}}

\newcommand{\Rd}{\ensuremath{\R^{d+1}}}
\newcommand{\Rn}{\ensuremath{\R^{2n-1}}}

\newcommand{\URn}{U\times\R^{2n+1}}

\newcommand{\Ca}[1]{\ensuremath{\mathcal{#1}}}

\newcommand{\cE}{\Ca{E}}

\newcommand{\cL}{\ensuremath{\mathcal{L}}}

\newcommand{\cV}{\ensuremath{H}}

\newcommand{\psivdo}{$\Psi_{H}$DO}
\newcommand{\psivdos}{$\Psi_{H}$DO's}
 
\newcommand{\pvdo}{\ensuremath{\Psi_{\cV}}}

\newcommand{\psido}{$\Psi$DO}

\newcommand{\psidos}{$\Psi$DO's}

\newcommand{\im}{\op{im}}

\begin{document}

\title{Complex powers of the contact Laplacian and the Baum-Connes conjecture for ${\rm SU}(n,1)$.}

\author{Rapha\"el Ponge}

\address{Max Planck Institute for Mathematics, Bonn, Germany}

\email{raphaelp@mpim-bonn.mpg.de}

\begin{abstract}
  This paper is an extended version of~\cite{Po:Crelle1b} where we point out and remedy a gap in the proof by Julg-Kasparov~\cite{JK:OKTGSU} of 
  the Baum-Connes conjecture for discrete subgroups of 
  $SU(n,1)$. In particular, here we explain in details why the non-microlocality of the Heisenberg calculus prevents us from implementing into 
  this framework the classical approach of Seeley to pseudodifferential complex powers, which was the main issue at stake in~\cite{Po:Crelle1b}.
\end{abstract}
   
 \maketitle  

\section{Introduction} 
A locally compact group $\Gamma$ satisfies the Baum-Connes conjecture when the Baum-Connes assembly map $\mu:K_{i}^{\op{top}}(\Gamma)\rightarrow 
K_{i}(C^{*}_{r}(\Gamma))$ is an isomorphism.  Here $K_{i}^{\op{top}}$ denotes the geometric $K$-group of Baum-Connes (when $\Gamma$ is torsion-free 
this the $K$-homology group 
$K_{i}(B\Gamma)$ of the classifiant space $B\Gamma$) and $K_{i}(C^{*}_{r}(\Gamma))$ denotes the analytic $K$-group of the reduced $C^{*}$-algebra of 
$\Gamma$. A stronger conjecture states that this holds with coefficients in any $C^{*}$-algebra acted on by 
$\Gamma$. Since its statement in the early 80's the Baum-Connes conjecture, with or without coefficients, has been shown for a variety of groups 
(see, e.g., \cite{HK:ETKTGWAPIHS}, \cite{Ju:CBCGSpn1}, \cite{Ka:LGKURCP}, \cite{La:KBABCBC}, \cite{MY:BCCHG}). 

In this paper we are concerned with the proof by Julg-Kasparov~\cite{JK:OKTGSU} of this conjecture for discrete subgroups of the complex Lorentz 
group $\op{SU}(n,1)$. The proof of Julg and Kasparov can be briefly summarized as follows. 

First, the proof can be reduced to showing that Kasparov's element $\gamma_{G}$ is equal to $1$ 
in the representation ring ${R}(G)$, that is, if $K$ is 
a maximal compact group of $G$ the restriction map ${R}(G)\rightarrow {R}(K)$ is an isomorphism. Second, the symmetric space $G/K$ is a complex 
hyperbolic space and under the Siegel map it is biholomorphic to the unit ball $B^{2n}\subset \C^{n}$ and its visual boundary is CR diffeomorphic to 
the unit sphere $S^{2n-1}$ equipped its standard CR structure. Julg and Kasparov further showed that $R(K)$ can be geometrically realized as 
${KK}_{G}(C(\overline{B}^{2n}),\C)$, where $C(\overline{B}^{2n})$ denotes the $C^{*}$-algebra of continuous functions on the closed unit ball 
$\overline{B}^{2n}$ and ${KK}_{G}$ is the equivariant ${KK}$ functor of Kasparov. Then they built a Fredholm 
module representing an element $\delta$ in ${KK}_{G}(C^{0}(\overline{B}^{2n}),\C)$ which is mapped to $\gamma$ in $KK_{G}(\C, \C)={R}(G)$ 
under the morphism induced by the map $\overline{B}^{2n}\rightarrow \{\text{pt}\}$. 

The construction of the element $\delta$ in ${KK}_{G}(C(\overline{B}^{2n}),\C)$ involves in a crucial manner  
the contact complex of Rumin~\cite{Ru:FDVC} on the unit sphere $S^{2n-1}$ endowed with its standard 
contact structure. In the contact setting the main geometric operators are not elliptic and the relevant pseudodifferential calculus to deal with them is 
the Heisenberg calculus of Beals-Greiner~\cite{BG:CHM} and Taylor~\cite{Ta:NCMA}. Then for constructing the Fredholm module representing 
$\delta$ Julg and Kasparov proved that the complex powers of the contact Laplacian are pseudodifferential operators in the Heisenberg calculus 
(see~\cite[Thm.~5.27]{JK:OKTGSU}). 

In~\cite{Se:CPEO} Seeley settled a general procedure for obtaining complex powers of elliptic operators as pseudodifferential operators. Its approach 
relied on constructing an asymptotic resolvent in a suitable class of classical \psidos\ with parameter. Similarly, Julg and Kasparov constructed an asymptotic 
resolvent in a class of Heisenberg \psidos\ with parameter (see~\cite[Thm.~5.25]{JK:OKTGSU}). We point out here that is class is not enough to allow 
us to carry out the rest of Seeley's arguments. In fact, we further show that we cannot carry out at all Seeley's approach in the setting of the Heisenberg 
calculus. Therefore, we have to rely on another approach to deal with the complex powers of the contact Laplacian. 
 
 This note is organized as follows. In Section~\ref{sec:Seeley-approach} we recall the outline of Seeley's approach to complex powers of elliptic 
 operators. In Section~\ref{sec:contact-Laplacian} we recall the construction of the contact complex and the associated contact Laplacian.
 In Section~\ref{sec:Heisenberg} we recall the definition of the Heisenberg calculus and stress out its lack of microlocality. In 
 Section~\ref{sec:Gap}  we explain the gap 
 in the proof of Theorem~5.27 of~\cite{JK:OKTGSU} and why the non-microlocality of the Heisenberg calculus prevents us from extending Seeley's approach
 to the setting of the  Heisenberg calculus. 
 
 \section{Seeley's approach to complex powers of elliptic operators}\label{sec:Seeley-approach}
 In this section we recall the main idea in the approach of Seeley~\cite{Se:CPEO} to complex powers of 
 elliptic operators  (see also~\cite{GS:WPPOAPSBP}, \cite{Sh:POST}). To simplify the exposition we let $M^{n}$ denote a compact manifold 
 equipped with a smooth density $>0$ and we let $P:C^{\infty}(M)\rightarrow C^{\infty}(M)$ be a positive elliptic differential operator of order $m$ with principal 
 symbol 
 $p_{m}(x,\xi)>0$.  Then for $\Re s<0$ we can write: 
 \begin{gather}
     P^{s}=\frac{i}{2\pi} \int_{\Gamma_{r}} \lambda^{s}(P-\lambda)^{-1}d\lambda, 
     \label{eq:SA.complex-powers-definition}\\
\Gamma_{r}=\{ \rho e^{i\pi}; \infty <\rho\leq r\}\cup\{ r e^{it}; 
\theta\geq t\geq \theta-2\pi \}\cup\{ \rho e^{-i\pi};  r\leq \rho\leq \infty\},
\label{eq:SA.contour}
 \end{gather}
 where $r$ is any real $<\lambda_{1}(P)$ and $\lambda_{1}(P)$ is the smallest non-zero eigenvalue of $P$. 
 
 To show that the formula above defines a \psido, Seeley constructs an asymptotic resolvent $Q(\lambda)$ as parametrix for $P-\lambda$ in a suitable \psido\ calculus 
 with parameter. More precisely, let $\Lambda \subset \C\setminus 0$ be an open angular sector $\theta <\arg \lambda<\theta'$ with $0<\theta<\pi<\theta'<2\pi$. 
 In the sequel we will say that a subset $\Theta \subset 
 [\R^{n}\times \C]\setminus 0$ is \emph{conic} when for any $t>0$ and any $ (\xi,\lambda)\in \Theta$ we have $(t\xi,t^{m}\lambda)\in \Theta$.
 For instance the subset $\R^{n}\times \Lambda \subset [\R^{n}\times \C]\setminus 0$ is conic. 
 
 Let $U \subset \R^{n}$ be a local chart for $M$. Then in $U$ the asymptotic resolvent has a symbol of the form
$q(x,\xi;\lambda) \sim \sum_{j\geq 0}q_{-m-j}(x,\xi;\lambda)$,
 where $\sim$ is taken in a suitable sense (see~\cite{Se:CPEO}) and there exists an open conic subset $\Theta\subset [\R^{n}\times \C]\setminus 0$ containing 
 $\R^{n}\times \Lambda$ such that each symbol $q_{-m-j}(x,\xi;\lambda)$ is smooth on $U\times \Theta$ and satisfies
 \begin{equation}
     q_{-m-j}(x,t\xi;t^{m}\lambda)=t^{-m-j}q(x,\xi;\lambda) \quad \forall t>0.
      \label{eq:SA.homogeneity}
 \end{equation}
 
 If $p(x,\xi)=\sum_{j=0}^{m}p_{m-j}(x,\xi)$ denotes the symbol of $P$ in the local chart $U$ then the symbol $ q(x,\xi;\lambda)$ satisfies
 $1\sim  (p(x,\xi)-\lambda)q(x,\xi;\lambda)+\sum_{\alpha \neq 0}\frac{1}{\alpha !}\partial_{\xi}^{\alpha}p(x,\xi)D_{x}^{\alpha}q(x,\xi;\lambda)$,
 from which we obtain
 \begin{gather}
     q_{-m}(x,\xi;\lambda)= (p_{m}(x,\xi)-\lambda)^{-1},
     \label{eq:SA.symbol-resolvent1}\\
      q_{-m-j}(x,\xi;\lambda)=-(p_{m}(x,\xi)-\lambda)^{-1}\!\! \sum_{\substack{|\alpha|+k+l=j,\\ l\neq j}} \!\! 
      \frac{1}{\alpha !}\partial_{\xi}^{\alpha}p_{m-k}(x,\xi)D_{x}^{\alpha}q_{-m-l}(x,\xi;\lambda).
     \label{eq:SA.symbol-resolvent2}
 \end{gather}
 
Set $\rho=\inf_{x \in U}\inf_{|\xi|=1}p(x,\xi)$.  Possibly by shrinking $U$ we may assume $\rho>0$. Let 
$    \Theta =[\R^{n}\times \Lambda] \cup\{(\xi;\lambda)\in \R^{n} \times \C; 0\leq |\lambda|<\rho|\xi|^{m}\}$.
 Then the formulas~(\ref{eq:SA.symbol-resolvent1}) and (\ref{eq:SA.symbol-resolvent2}) show that each symbol $q_{-m-j}(x,\xi;\lambda)$ is well defined and 
 smooth on $U\times \Theta$ and is homogeneous in 
 the sense of~(\ref{eq:SA.homogeneity}). Furthermore, it is analytic with respect to $\lambda$. Therefore, for $\Re s<0$ we define a smooth function on $U\times 
( \R^{n}\setminus 0)$ by letting 
\begin{equation}
    p_{s,ms-j}(x,\xi)= \frac{i}{2\pi} \int_{\Gamma(\xi)} \lambda^{s}q_{-m-j}(x,\xi;\lambda)d\lambda,
\end{equation}
where $\Gamma(\xi)$ is contour $\Gamma_{r}$ in~(\ref{eq:SA.contour}) with $r=\frac{1}{2}\rho|\xi|^{m}$. Moreover, one can check that 
$p_{s,ms-j}(x,t\xi)=t^{ms-j}p_{s,ms-j}(x,\xi)$ for any $t>0$, i.e., $p_{s,ms-j}(x,\xi)$ is a homogeneous symbol of degree $ms-j$. 

It can also be shown that on the chart $U$ the operators $P^{s}$ is a \psido\ with symbol $p_{s} \sim \sum_{j \geq 0}p_{s,ms-j}(x,\xi)$.  This is true on any 
local chart so, as one can check 
that the Schwartz kernel of $P^{s}$ is smooth off the diagonal of $M\times M$, it follows that $P^{s}$ is a \psido\ of order $ms$ for $\Re s<0$. 

Furthermore, for $k=1,2,\ldots$ and $\Re s<k$ we have $P^{s}=P^{k}P^{s-k}$. Here $P^{k}$ is a differential operator of order $k$ and as we have $\Re s -k<0$ we 
know that $P^{s-k}$ is a \psido\ of order $m(s-k)$. Therefore, $P^{s}$ is a \psido\ of order $ms$ for $\Re s<k$ and any $k=1,2,,\ldots$, i.e, 
$P^{s}$ is a \psido\ of order $ms$ for any $s\in \C$. 

\section{Rumin's contact complex}\label{sec:contact-Laplacian}
 Let $(M^{2n-1},H)$ be an orientable contact manifold, so that there exists a global nonvanishing contact form $\theta$ such that $H=\ker \theta$ and 
 $d\theta_{|_{H}}$ is nondegenerate. Let $X_{0}$ be the 
Reeb vector field of $\theta$, so that $\iota_{X_{0}}\theta=1$ and $\iota_{X_{0}}d\theta=0$. In addition, we let $J$ be a calibrated almost complex structure on 
$H$, so that we have $\delta\theta(X,JX)>0$ for any non-zero section of $H$, and  we endow $TM$ with the Riemannian metric 
$g_{\theta,J}=d\theta(.,J.)+\theta^{2}$.

Observe that the splitting $TM=H\oplus \R X_{0}$ allows us to identify  
$H^{*}$ with the annihilator of $X_{0}$ in $T^{*}M$. More generally, by identifying $\Lambda^{k}_{\C}H^{*}$ with $\ker \iota_{X_{0}}$ we get the 
splitting,
\begin{equation}
    \Lambda^{*}_{\C}TM=(\bigoplus_{k=0}^{2n}\Lambda^{k}_{\C}H^{*}) \oplus (\bigoplus_{k=0}^{2n} \theta\wedge \Lambda^{k}_{\C}H^{*}).
     \label{eq:contact.decomposition-forms}
\end{equation}
 
For any horizontal form $\eta\in C^{\infty}(M,\Lambda^{k}_{\C}H^{*})$ we can write
$d\eta= d_{b}\eta+\theta \wedge \cL_{X_{0}}\eta$,
where $d_{b}\eta$ is the component of $d\eta$ in $\Lambda^{k}_{\C}H^{*}$. This does not provide us with a complex, for we have 
$d_{b}^{2}=-\cL_{X_{0}}\varepsilon(d\theta)=-\varepsilon (d\theta)\cL_{X_{0}}$, where $\varepsilon(d\theta)$ denotes the exterior multiplication 
by $d\theta$.

The contact complex of Rumin~\cite{Ru:FDVC} is an attempt to get 
a complex of horizontal differential forms by forcing the equalities $d_{b}^{2}=0$ and $(d^{*}_{b})^{2}=0$.

A natural way to modify $d_{b}$ to get the equality $d_{b}^{2}=0$ is to restrict 
$d_{b}$ to the subbundle $\Lambda^{*}_{2}:=\ker \varepsilon(d\theta) \cap \Lambda^{*}_{\C}H^{*}$, since the latter
is closed under $d_{b}$ and is annihilated by $d_{b}^{2}$. 

Similarly, we get the equality $(d_{b}^{*})^{2}=0$ by restricting $d^{*}_{b}$ to the subbundle 
$\Lambda^{*}_{1}:=\ker \iota(d\theta)\cap \Lambda^{*}_{\C}H^{*}=(\im \varepsilon(d\theta))^{\perp}\cap \Lambda^{*}_{\C}H^{*}$, where 
$\iota(d\theta)$ denotes the interior product 
with $d\theta$. This amounts to replace $d_{b}$ by $\pi_{1}\circ d_{b}$, where $\pi_{1}$ is the orthogonal projection onto $\Lambda^{*}_{1}$.

In fact, since $d\theta$ is nondegenerate on $H$ the operator $\varepsilon(d\theta):\Lambda^{k}_{\C}H^{*}\rightarrow \Lambda^{k+2}_{\C}H^{*}$  is 
injective for $k\leq n-1$ and surjective for $k\geq n+1$. This implies that $\Lambda_{2}^{k}=0$ for $k\leq n$ and $\Lambda_{1}^{k}=0$ for $k\geq n+1$. 
Therefore, we only have two halves of complexes. 

As observed by Rumin~\cite{Ru:FDVC} we get a full complex by connecting 
the two halves by means of the operator $D_{R,n}:C^{\infty}(M,\Lambda_{\C}^{n}H^{*}) \rightarrow 
C^{\infty}(M,\Lambda_{\C}^{n}H^{*})$ such that 
\begin{equation}
    D_{R,n}=\cL_{X_{0}}+d_{b,n-1}\varepsilon(d\theta)^{-1}d_{b,n},
\end{equation}
where $\varepsilon(d\theta)^{-1}$ is the inverse of $\varepsilon(d\theta):\Lambda^{n-1}_{\C}H^{*}\rightarrow \Lambda^{n+1}_{\C}H^{*}$. Notice that 
$D_{R,n}$ is second order differential 
operator.  This allows us to get the contact complex, 
\begin{equation}
    C^{\infty}(M)\stackrel{d_{R,0}}{\rightarrow}
    \ldots 
    C^{\infty}(M,\Lambda^{n})\stackrel{D_{R,n}}{\rightarrow} C^{\infty}(M,\Lambda^{n}) 
    \ldots \stackrel{d_{R,2n-1}}{\rightarrow} C^{\infty}(M,\Lambda^{2n}).
     \label{eq:contact-complex}
\end{equation}
where $d_{R,k}$ agrees with $\pi_{1}\circ d_{b}$ for $k=0,\ldots,n-1$ and  with $d_{R,k}=d_{b}$ otherwise. 

The contact Laplacian is defined as follows. In degree $k\neq n$ this is the differential operator 
$\Delta_{R,k}:C^{\infty}(M,\Lambda^{k})\rightarrow C^{\infty}(M,\Lambda^{k})$ such that
\begin{equation}
    \Delta_{R,k}=\left\{
    \begin{array}{ll}
        (n-k)d_{R,k-1}d^{*}_{R,k}+(n-k+1) d^{*}_{R,k+1}d_{R,k},& \text{$k=0,\ldots,n-1$},\\
         (k-n-1)d_{R,k-1}d^{*}_{R,k}+(k-n) d^{*}_{R,k+1}d_{R,k},& \text{$k=n+1,\ldots,2n$}.
         \label{eq:contact-Laplacian1}
    \end{array}\right.
\end{equation}
For $k=n$ we have the differential operators $\Delta_{R,nj}:C^{\infty}(M,\Lambda_{j}^{n})\rightarrow C^{\infty}(M,\Lambda^{n}_{j})$, $j=1,2$, 
given by the formulas, 
\begin{equation}
    \Delta_{R,n1}= (d_{R,n-1}d^{*}_{R,n})^{2}+D_{R,n}^{*}D_{R,n}, \quad   \Delta_{R,n2}=D_{R,n}D_{R,n}^{*}+  (d^{*}_{R,n+1}d_{R,n}).
         \label{eq:contact-Laplacian2}
\end{equation}

Observe that $\Delta_{R,k}$, $k\neq n$, is a differential operator order $2$, whereas $\Delta_{Rn1}$ and $\Delta_{Rn2}$ are differential operators of 
order $4$. Moreover, Rumin~\cite{Ru:FDVC} proved that in every degree the contact Laplacian is maximal hypoelliptic. 

\section{Heisenberg calculus}\label{sec:Heisenberg}
Let $(M^{2n-1},H)$ be a contact manifold. In this setting the natural operators like the contact Laplacian are not elliptic, so the standard 
 pseudodifferential calculus does not apply. The substitute is provided by the Heisenberg calculus of Beals-Greiner~\cite{BG:CHM} and Taylor~\cite{Ta:NCMA}. 
 The idea, which goes back to Elias 
 Stein, is to construct a class of pseudodifferential operators, the \psivdos, which at each point are modelled by left-convolutions operators on the Heisenberg 
 group $\bH^{2n-1}$. This motivated by the fact that in a suitable sense $\bH^{2n-1}$ is tangent to a contact manifold at each of its points. 
 
  Locally the \psidos\ can be described as follows. Let $U \subset \Rn$ be a local chart together with a frame $X_{0},\ldots,X_{2n}$ such that 
  $X_{1},\ldots,X_{2n}$ span $H$. Such a chart is called a Heisenberg chart. Moreover, on $\Rn$ we consider the dilations, 
\begin{equation}
     t.\xi=(t^{2}\xi_{0},t\xi_{1},\ldots,t\xi_{2n}), \qquad \xi\in \Rd, \quad t>0. 
     \label{eq:HC.Heisenberg-dilations}
\end{equation}
 
\begin{definition}1)  $S_{m}(\URn)$, $m\in\C$, is the space of functions 
    $p(x,\xi)$ in $C^{\infty}(U\times(\Rn\setminus 0))$ such that $p(x,t.\xi)=t^m p(x,\xi)$ for any $t>0$.\smallskip

2) $S^m(\URn)$,  $m\in\C$, consists of functions  $p\in C^{\infty}(\URn)$ with
an asymptotic expansion $ p \sim \sum_{j\geq 0} p_{m-j}$, $p_{k}\in S_{k}(\URn)$, in the sense that, for any integer $N$ and 
for any compact $K \subset U$, we have
\begin{equation}
    | \partial^\alpha_{x}\partial^\beta_{\xi}(p-\sum_{j<N}p_{m-j})(x,\xi)| \leq 
    C_{\alpha\beta NK}\|\xi\|^{\Re m-\brak\beta -N}, \qquad  x\in K, \quad \|\xi \| \geq 1,
    \label{eq:NCRP.asymptotic-expansion-symbols}
\end{equation}
where we have let $\brak\beta=2\beta_{0}+\beta_{1}+\ldots+\beta_{2n}$ and 
 $\|\xi\|=(\xi_{0}^{2}+\xi_{1}^{4}+\ldots+\xi_{2n}^{4})^{1/4}$.
\end{definition}

Next, for $j=0,\ldots,2n$ let  $\sigma_{j}(x,\xi)$ denote the symbol (in the 
classical sense) of the vector field $\frac{1}{i}X_{j}$  and set  $\sigma=(\sigma_{0},\ldots,\sigma_{2n})$. Then for $p \in S^{m}(\URn)$ we let $p(x,-iX)$ be the 
continuous linear operator from $C^{\infty}_{c}(U)$ to $C^{\infty}(U)$ such that 
\begin{equation}
          p(x,-iX)f(x)= (2\pi)^{-(d+1)} \int e^{ix.\xi} p(x,\sigma(x,\xi))\hat{f}(\xi)d\xi,
    \qquad f\in C^{\infty}_{c}(U).
\end{equation}

\begin{definition}
   $\pvdo^{m}(U)$, $m\in \C$, consists of operators $P:C^{\infty}_{c}(U)\rightarrow C^{\infty}(U)$ which are of the form
$P= p(x,-iX)+R$ for some $p$ in $S^{m}(\URn)$, called the symbol of $P$, and some smoothing operator $R$.
\end{definition}
  
The class of \psivdos\ is invariant under changes of Heisenberg chart (see~\cite[Sect.~16]{BG:CHM}, \cite[Appendix A]{Po:MAMS1}), so we may 
 extend the definition of \psivdos\ to an arbitrary Heisenberg manifold $(M,H)$ and let them act on sections of a vector bundle $\cE$ over $M$. 
 We let $\pvdo^{m}(M,\cE)$ denote the class of \psivdos\ of order $m$ on $M$ acting on sections 
 of $\cE$. 

 In addition, the Heisenberg calculus possesses a full symbolic calculus which allows us to construct parametrices for hypoelliptic operators. 
 In the classical pseudodifferential calculus the relevant product at  the level of homogeneous symbols is the pointwise product of functions which, 
 under the Fourier transform, corresponds to the convolution on the Abelian group $\R^{2n-1}$. 
 
 Similarly, for homogeneous Heisenberg symbols the relevant product comes from the convolution on the Heisenberg group $\bH^{2n-1}$. More precisely, 
 let $U\subset \Rn$ be a Heisenberg chart  as above and for $m \in \C$ let $S_{m}(\Rn)$ denotes the closed subspace of $C^{\infty}(\Rn\setminus 0)$ 
 consisting in functions $p(\xi)$ such that $p(t.\xi)=t^{m}p(\xi)$ for any $t>0$. Then for each $x\in U$ we get a bilinear product
 \begin{equation}
     *^{x}:S_{m_{1}}(\Rn)\times S_{m_{2}}(\Rn) \longrightarrow S_{m_{1}+m_{2}}(\Rn).
 \end{equation}
 This product depends smoothly on $x$, so we may define the bilinear product, 
\begin{gather}
    *:S_{m_{1}}(\URn)\times S_{m_{2}}(\URn) \longrightarrow S_{m_{1}+m_{2}}(\URn),
    \label{eq:HC.composition-homogeneous-symbols1}\\
    p_{m_{1}}*p_{m_{2}}(x,\xi)=[p_{m_{1}}(x,.)*^{x}p_{m_{2}}(x,.)](\xi), \qquad p_{m_{j}} \in S_{m_{j}}(\URn).
     \label{eq:HC.composition-homogeneous-symbols2}
\end{gather}
Then we have: 
\begin{proposition}[{\cite[Thm.~14.7]{BG:CHM}}] \label{prop:HC.composition}
    For $j=1,2$ let $P_{j}\in \pvdo^{m_{j}}(U)$ have  symbol $p_{j}\sim \sum_{k\geq 0} p_{j,m_{j}-k}$ and assume  
    that one of these operators is properly supported. Then the operator  
$P=P_{1}P_{2}$ is a \psivdo\ of order $m_{1}+m_{2}$ and has symbol  $p\sim \sum_{k\geq 0} p_{m_{1}+m_{2}-k}$, with  
\begin{equation}
     p_{m_{1}+m_{2}-k} = \sum_{k_{1}+k_{2}\leq k} \sum_{\alpha,\beta,\gamma,\delta}^{(k-k_{1}-k_{2})}
            h_{\alpha\beta\gamma\delta}  (D_{\xi}^\delta p_{1,m_{1}-k_{1}})* (\xi^\gamma 
            \partial_{x}^\alpha \partial_{\xi}^\beta p_{2,m_{2}-k_{2}}), 
            \label{eq:HC.composition}
\end{equation}
where $\underset{\alpha\beta\gamma\delta}{\overset{\scriptstyle{(l)}}{\sum}}$ denotes the sum over all the indices such that 
$|\alpha|+|\beta| \leq \brak\beta -\brak\gamma+\brak\delta = l$ and $|\beta|=|\gamma|$, and the functions 
$h_{\alpha\beta\gamma\delta}(x)$'s are  polynomials in the derivatives of the coefficients of 
the vector fields $X_{0},\ldots,X_{d}$.
\end{proposition}

This result allows us to carry out the classical parametrix construction, provided we can invert the principal symbol with respect to the product 
$*$. This may be difficult in practice because this is not anymore the pointwise product of functions, but 
this can be completely determined in terms of a representation theoretic criterion, 
the so-called Rockland condition (see~\cite{Po:MAMS1}). For instance, in every degree the contact Laplacian satisfies this condition and 
admits a parametrix in the Heisenberg calculus (see~\cite{JK:OKTGSU},~\cite[Sect.~3.5]{Po:MAMS1}).

Finally, it should be stressed out that in~(\ref{eq:HC.composition-homogeneous-symbols1})--(\ref{eq:HC.composition-homogeneous-symbols2}) 
the $x$-values of $p_{m_{1}}*p_{m_{2}}$ depend only on that of $p_{m_{1}}$ and $p_{m_{2}}$, i.e., 
$p_{m_{1}}*p_{m_{2}}(x,.)$ depends only on $x$ and $p_{m_{1}}(x,.)$ and $p_{m_{2}}(x,.)$. However, for any $(x,\xi) \in U\times (\Rn\setminus 0)$ the 
value of $p_{m_{1}}*p_{m_{2}}(x,\xi)$ depends on all the values of $p_{m_{1}}(x,\xi')$ and $p_{m_{2}}(x,\xi')$ as $\xi'$ ranges over $\Rn\setminus 
0$. Thus the Heisenberg calculus is local, but is not microlocal.

\section{Complex powers of the contact Laplacian}\label{sec:Gap}
Let $S^{2n-1}\subset \C^{n}$ be the unit sphere equipped with its standard contact structure $H=\ker \theta$ with 
$\theta=\sum_{j=0}^{n} \overline{z}_{j}dz_{j}$. Since $H$ is invariant under the multiplication by $i$ the 
complex structure of $\C^{n+1}$ induces a complex structure $J$ on $H$. This is a calibrated complex structure on $H$ and 
we then endow $TM$ with the 
Riemannian metric $g_{\theta,J}=d\theta(.,J.)+\theta^{2}$.

Let $\Delta_{R}$ be the contact Laplacian of $S^{2n-1}$. In the proof of the Baum-Connes conjecture with coefficients for discrete subgroups of 
$SU(n,1)$ in~\cite{JK:OKTGSU} the following 
result is needed. 

\begin{theorem}[{\cite[Thm.~5.27]{JK:OKTGSU}}]\label{thm:CPCL.Thm}
    Let $s \in \C$. Then:\smallskip 
    
    1) The operator $\Delta_{R,k}^{s}$, $k \neq n$,  is \psivdo\ of order $2s$;\smallskip
    
    2) The operator $\Delta_{R,nj}^{s}$, $j=1,2$, is a \psivdo\ of order $4s$.
\end{theorem}

We would like to point out a gap in the proof of this result in~\cite{JK:OKTGSU}. For simplicity we will explain it in degree $k \neq n$, but the argument works  in 
degree $k=n$ as well. The idea in the proof in~\cite{JK:OKTGSU} is to carry out Seeley's approach in the Heisenberg setting. 
As in~(\ref{eq:SA.complex-powers-definition}) for $\Re s<0$ we have 
\begin{equation}
    \Delta_{R,k}^{s}=\frac{i}{2\pi} \int_{\Gamma} \lambda^{s}(\Delta_{R,k}-\lambda)^{-1}d\lambda, 
     \label{eq:Gap.complex-power}
\end{equation}
where $\Gamma$ is as in~(\ref{eq:SA.contour}). Let $\Lambda \subset \C\setminus 0$ be an open angular sector $\theta <\arg \lambda<\theta'$ 
with $0<\theta<\pi<\theta'<2\pi$. Then Julg and Kasparov 
showed that the resolvent $(\Delta_{R,k}-\lambda)^{-1}$ belongs to a class $\pvdo^{-2}(S^{2n-1},\Lambda^{k};\Lambda)$ of \psivdos\ parametrized by 
$\Lambda$ (see~\cite[Thm.~5.25]{JK:OKTGSU}). 
In particular, in a Heisenberg chart $U \subset \Rn$ the resolvent $(\Delta_{R,k}-\lambda)^{-1}$ has a 
symbol $q(x,\xi;\lambda)$ in $C^{\infty}(\URn \times \Lambda)$ with an expansion
$q(x,\xi;\lambda) \sim \sum_{j\geq 0} q_{-2-j}(x,\xi;\lambda)$,
where $\sim$ is taken in a suitable sense and $q_{-2-j}(x,\xi;\lambda)\in C^{\infty}(U\times (\Rn \setminus 0)\times \Lambda))$ is such that for any 
$t>0$ we have 
$q_{-2-j}(x,t.\xi;t^{2}\lambda)=t^{-2-j}q_{-2-j}(x,\xi;\lambda)$. 

Notice that  the symbol $q(x,\xi;\lambda)$ does not make sense for $\lambda>0$, because only angular sectors contained in $\C \setminus [0,\infty)$ 
are considered in~\cite[Thm.~5.25]{JK:OKTGSU}. However, we have to allow the homogeneous components of the symbol of $(\Delta_{R,k}-\lambda)^{-1}$ 
to be defined for 
some $(x,\xi;\lambda)\in U\times (\Rn\times0)\times \C $ with $\lambda>0$ because the contour $\Gamma$ in~(\ref{eq:Gap.complex-power}) 
always crosses the positive real axis. 

This leads us to a discrepancy in~\cite[p.~128, line 9 from bottom]{JK:OKTGSU} when the authors claim that their Theorem~5.25 insures us that 
the resolvent $(\Delta_{R,k}-\lambda)^{-1}$ belongs to a class $\pvdo^{-2}(S^{2n-1},\Lambda^{k};\Lambda)$ for some subset $\Lambda$ containing the contour 
$\Gamma$. In particular, the authors cannot claim that, in a Heisenberg chart $U \subset \Rn$, the symbol of $\Delta_{R,k}^{s}$ is given by the formula, 
\begin{equation}
    p_{s}(x,\xi)=  \frac{i}{2\pi} \int_{\Gamma} \lambda^{s}q(x,\xi;\lambda)d\lambda,
\end{equation}
where $q(x,\xi;\lambda)$ denotes the symbol of $(\Delta_{R,k}-\lambda)^{-1}$, since the homogeneous components of the latter symbol don't make sense 
for $\lambda>0$. 

In fact, we cannot implement at all Seeley's approach to the setting of the Heisenberg calculus. More precisely, 
to carry out Seeley's approach we have to construct an asymptotic resolvent for $\Delta_{R,k}$ in a class of \psivdos\ with 
parameter associated to an angular sector $\Lambda$ as above and given in a local Heisenberg chart $U \subset \Rn$ by parametric symbols, 
$q(x,\xi;\lambda) \sim \sum_{j\geq 0}q_{-m-j}(x,\xi;\lambda)$, 
where $\sim$ is taken in a suitable sense and there exists an open conic subset $\Theta\subset [\R^{n}\times \C]\setminus 0$ containing 
 $\R^{n}\times \Lambda$ such that each symbol $q_{-m-j}(x,\xi;\lambda)$ is smooth on $U\times \Theta$ and satisfies 
 $q_{-m-j}(x,t.\xi;t^{2}\lambda)=t^{-2-j}q(x,\xi;\lambda)$ for any $t>0$. 
 
 If we let $p(x,\xi)=\sum p_{2-j}(x,\xi)$ be the symbol of $\Delta_{R,k}$ in the Heisenberg chart, 
 then by Proposition~\ref{prop:HC.composition} we have
  \begin{equation}
     1 \sim \sum_{j \geq 0}  \sum_{k+l\leq j} \sum_{\alpha,\beta,\gamma,\delta}^{(j-k-l)}
             h_{\alpha\beta\gamma\delta}(x)  (D_{\xi}^\delta p_{2-k})* (\xi^\gamma 
             \partial_{x}^\alpha \partial_{\xi}^\beta q_{-2-l})(x,\xi;\lambda),
  \end{equation}
where the notation is the same as in~(\ref{eq:HC.composition}). Therefore, we get
\begin{gather}
    q_{-2}(x,\xi;\lambda)=(p_{2}-\lambda)^{*-1}(x,\xi;\lambda)
    \label{eq:Gap.symbol-resolvent1}  \\
    q_{-2-j}(x,\xi;\lambda)=-\sum_{\substack{k+l\leq j,\\ l\neq j}} \sum_{\alpha,\beta,\gamma,\delta}^{(j-k-l)}
            h_{\alpha\beta\gamma\delta}(x)  q_{-2}*(D_{\xi}^\delta p_{2-k})* (\xi^\gamma 
            \partial_{x}^\alpha \partial_{\xi}^\beta q_{-2-l})(x,\xi;\lambda),
    \label{eq:Gap.symbol-resolvent2} 
\end{gather}
where $(p_{2}-\lambda)^{*-1}$ denotes the inverse of $p_{2}-\lambda$ with respect to the product $*$. 

If $q_{1}(x,\xi;\lambda)$ and $q_{2}(x,\xi;\lambda)$ are two homogeneous Heisenberg symbols with parameter then the product $q_{1}$ and $q_{2}$ 
should be defined as 
\begin{equation}
    q_{1}*q_{2}(x,\xi;\lambda)=[q_{1}(x,.;\lambda)*^{x}q_{2}(x,.;\lambda)](\xi).
\end{equation}
As mentioned at the end of Section~\ref{sec:Heisenberg} the definition of $[q_{1}(x,.;\lambda)*^{x}q_{2}(x,.;\lambda)](\xi)$ depends on all the values of 
$q_{1}(x,\xi';\lambda)$ and $q_{2}(x,\xi';\lambda)$ as $\xi'$ ranges over $\Rn\setminus 0$. For a parameter 
$\lambda>0$ the symbols $q_{1}(x,\xi;\lambda)$ and 
$q_{2}(x,\xi;\lambda)$ are only defined for $\xi$ in $\{\xi; (x,\xi;\lambda)\in \Theta\}$ which does not agree with $\Rn\setminus 0$, so we cannot 
define $q_{1}*q_{2}(x,\xi;\lambda)$ for $\lambda>0$. Therefore, the formula~(\ref{eq:Gap.symbol-resolvent2}) does not make sense for $\lambda>0$. 

All this shows that the non-microlocality of the Heisenberg calculus prevents us from implementing Seeley's approach into the setting of the Heisenberg 
calculus. 

Therefore, we have to rely on other approaches to deal with complex powers within the framework of this calculus. This is done for instance 
 in~\cite[Sect.~5.3]{Po:MAMS1} for a large class of operators, including the contact Laplacian on any compact contact manifold. This
 allows us to complete the proof in~\cite{JK:OKTGSU} of the Baum-Connes conjecture with coefficients for discrete subgroups of ${SU}(n,1)$.
 
 The approach in~\cite{Po:MAMS1} is based on combining Mellin's formula for the complex powers with the pseudodifferential representation of the 
 heat kernel of~\cite{BGS:HECRM}. Furthermore, the above-mentioned issues with considering Heisenberg \psidos\ with parameter associated to 
 homogeneous symbols with parameter disappear when we substitute the latter by \emph{almost} homogeneous symbols symbols with parameter 
 (see~\cite{Po:PhD}, \cite{Po:CPDE1}). This provides us with an alternative approach to deal with complex power within the framework of the Heisenberg 
 calculus. 

It is also possible to construct complex powers in this setting of the Heisenberg by means of a pseudodifferential representation of the resolvent. 
Instead of using homogeneous symbols with parameters as above,  we can use \emph{almost} homogeneous symbols with parameter. This is more suitable 
for domains containing contour $\Gamma$ as in~(\ref{eq:SA.contour}) and~(\ref{eq:Gap.complex-power}) and allows us to resolve all the issues 
raised above. This is done for sublaplacians in~\cite{Po:PhD} and for more general operators in~\cite{Po:CPDE1}.

\end{document}